\documentclass[11pt]{article}
\usepackage [a4paper,left=2.5cm,bottom=2.5cm,right=2.5cm,top=2.5cm]{geometry}
\usepackage[english]{babel}
\usepackage{graphicx,color,subfigure,multirow, here,ulem}
\usepackage{mathtools}
\usepackage[utf8]{inputenc} 
\usepackage[T1]{fontenc}
\usepackage[english]{babel} 
\usepackage[numbers]{natbib}
\usepackage{amsmath, amsfonts, amsthm,amssymb,here,dsfont, hyperref, enumitem}
\usepackage{natbib}
\newcommand{\argmin}[1]{\underset{#1}{\operatorname{arg}\!\operatorname{min}}\;}
\newcommand{\one}{\mathds{1}}

\newcommand{\w}{\widehat}
\newcommand{\wt}{\widetilde}

\newcommand{\pen}{\text{pen}}
\makeatletter

\@addtoreset{equation}{section}
\makeatother

\definecolor{cadmiumgreen}{rgb}{0.0, 0.50, 0.24}

\begin{document}

\title{Neuronal Network Inference and Membrane Potential Model using Multivariate Hawkes Processes}
%{Reconstructing the membrane potential of a neuron surrounded by a neural network}

%\begin{center}
%\small 
\author{
Anna Bonnet$^{(1)}$, Charlotte Dion-Blanc$^{(1)}$\footnote{Corresponding author: \url{charlotte.dion_blanc@sorbonne-universite.fr}}, François Gindraud$^{(2)}$ and Sarah Lemler$^{(3)}$ \\
$^{(1)}$ Sorbonne Université, UMR CNRS 8001, LPSM, 75005 Paris, France\\
$^{(2)}$ Université Lyon 1, CNRS, LBBE UMR 5558, F-69622 Villeurbanne, France\\
$^{(3)}$ Laboratoire MICS, \'Ecole CentraleSup\'elec, Université Paris-Saclay
}
%\end{center}

\maketitle
  \begin{abstract}
In this work, we propose to catch the complexity of the membrane potential's dynamic of a motoneuron between its spikes, taking into account the spikes from other neurons around. Our approach relies on two types of data: extracellular recordings of multiple spikes trains and intracellular recordings of the membrane potential of a central neuron. Our main contribution is to provide a unified framework and a complete pipeline to analyze neuronal activity from data extraction to statistical inference. The first step of the procedure is to select a subnetwork of neurons impacting the central neuron: we use a multivariate Hawkes process to model the spike trains of all neurons and compare two sparse inference procedures to identify the connectivity graph. Then we infer a jump-diffusion dynamic in which jumps are driven from a Hawkes process, the occurrences of which correspond to the spike trains of the aforementioned subset of neurons that interact with the central neuron. We validate the Hawkes model with a goodness-of-fit test and we show that taking into account the information from the connectivity graph improves the inference of the jump-diffusion process. The entire code has been developed and is freely available on GitHub.
   \end{abstract}



\section{Introduction}\label{sec:intro}

Motoneurons are the center of the movement commands. 
%In practice, two types of neuronal signals can be measured: an intracellular signal by introducing an electrode inside a neuron, the membrane potential of a neuron is then recovered; or an extracellular signal, using a multi-electrode, the instants (called spikes) at which information is sent by several neurons simultaneously are recovered, the signal is then processed to find out which neurons correspond to which spikes, this is the spike sorting. 
%\ab{\sout{ In this work, we investigate the evolution of the dynamic of the membrane potential of a fixed motoneuron from the spike trains of a neuron network around.}}
In practice, several types of signals can be recovered: in particular, an intracellular signal, obtained by introducing an electrode inside the neuron and which measures its membrane potential and also an extracellular signal, obtained with a multi-electrode, that contains the moments (called spikes) at which information is sent by several neurons simultaneously. This extracellular signal is then processed to determine the correspondence between neurons and spikes: this step is called the spike sorting. 

Usually, neuronal studies focus either on the intracellular signal to model the dynamics of the membrane potential of a neuron or on the extracellular signals to recover the graph of connectivity between several neurons. In this paper, we aim at taking advantage of both signals to infer first the interactions between neurons and then the dynamic of the membrane potential of a fixed neuron, accounting for these estimated interactions. %spike trains of a neuron network of the motor system and a membrane potential from one of them, assuming in the center of the network.
The data recorded on turtles have already been used in \cite{berg2013} and \cite{radosevic}, and are presented in detail in Section \ref{sec:data}. We have access to the spike trains of 249 neurons and the membrane potential of one neuron, which we call the central neuron. Each experiment has been conducted 10 times, so we have 10 trials of these signals. 
In order to increase the understanding of the link between the macro phenomenon and the micro behavior of the neurons, our goal is to answer the two following questions: 
\begin{itemize}
    \item What is the strength of the connections between the neurons?
    \item How does the network around a neuron of interest impact its membrane potential?
\end{itemize} 
First, the paper aims at providing a graph of connectivity representing
the rhythmic scratch-like network activity in the turtle. The second goal is to fit a dynamic model for the membrane potential of a neuron impacted by this interaction network between the spikes.

\subsection{Related works}

When dealing with neuron data, two standard centers of interest are either the interactions between neurons or the evolution of the dynamic of the membrane potential. 
Regarding the first aspect, neurons interact mainly through synapses.
When the neuron sends a message, it produces a spike, corresponding to a substantial increase of its membrane potential. The times when these spikes happen, called the spike trains, are then recorded by neurophysiologists to be analyzed.

There is extensive literature on the study of multiple spike trains. Some classical approaches are used to study the correlation between pairs of neurons, among them, we can mention, for example, the histogram-based method to recover the cross-correlation function \cite{PERKEL1967419} or the joint Peri-Stimulus Time Histogram (PSTH) \cite{gerstein1969simultaneously}.

More recently, methods based on point processes have been proposed to study a neuron network with more than two neurons. To cite a few, \cite{pouzat} investigates a model based on Poisson processes, \cite{rebesco} chooses a doubly stochastic Poisson process to infer functional connectivity. Nevertheless, if inhomogeneous Poisson processes are adapted to model time heterogeneity, they are not well suited to handle dependency between points.
%\\https://fr.overleaf.com/project/604a71ae91e6cf5963a00531
%doubly stochastic Poisson process \cite{rebesco}: Inferred Functional Connectivity (IFC) algorithm (bayesian algorithm) to estimate the way neurons interact with one another and obtain  the functional connectivity. The method incorporates constraints that penalize a large number of connections, the approach provides a parsimonious and robust description of network connectivity. Used to obtain connectivity changes because of stimulations.\\
Some papers do not impose a time series model but use Generalized Linear Model (GLM), see for example
 \cite{roudi} or  \cite{kobayashi2019reconstructing}. The last reference focuses on reconstructing connectivity circuit from GLM model on classical cross-correlogram (see \emph{e.g.} \cite{CCH}), for a large-scale network.
 
Hawkes processes have also been proposed to model jointly the spike trains of several neurons and then infer the interactions between these neurons (see \cite{RBRTM}, \cite{reynaud2014goodness}  and \cite{DGL2016}).
Such point processes are commonly used to model a complex phenomenon with strong time connections as earthquakes \citep{ogata88}, Twitter networks \cite{twit}, order books \cite{BACRYFINANCE} or spatial dependencies, for instance in genomics \cite{bonnet}.
Hawkes processes are historically characterized by an auto-excitatory behavior although they can model other types of dependencies, in particular inhibitive effects \cite{BMS21}.
%, which fits well with the neurons' behavior. 
%In neuroscience, the Hawkes model has been proposed to model jointly the spike trains of several neurons,  see for example \cite{RBRTM}, \cite{reynaud2014goodness}  and \cite{DGL2016}.
 In \cite{GLMhawkes},
 the authors propose a stability framework for the data-driven (univariate) non-linear Hawkes process investigated with GLM. 
 %: Point Process Generalized Linear Model, passes common goodness-of-fit tests, but is unstable and can lead to divergent firing rates. In the paper \cite{GLMhawkes}, they propose a stability framework for data-driven (univariate) Non-linear Hawkes PP-GLMs, a quasi-renewal approximation that can be used alongside other goodness-of-fit tests to study the stability and assess model adequacy. 
 %(In the perspectives) There exists an optimally L1-regularized solution that is predicted to be stable: therefore, strong L1-regularization might be an alternative approach to model stabilization.
Hawkes processes for spike trains are also used in \cite{lambert}, which is based on the estimator presented in \cite{Hansen15}. It presents a methodology to study spike trains for small networks, modeling the data with a multivariate Hawkes process. 
%In the following we consider this procedure on the data. 
%\sout{We also work with the algorithm of \cite{ADM4} implemented in Python library tick to simulate and study Hawkes processes. The inference procedure developed in \cite{ADM4} is based on the likelihood and is fitted for  large networks.}

%\sout{To validate the choice of the Hawkes process to model the interaction between neurons on our data, we also apply the test procedure introduced in \cite{reynaud2014goodness}. This is an asymptotic test based on the Kolmogrorov statistics, which is detailed in Section \ref{sec:test}.
%\textcolor{red}{à voir si on laisse ça ici :}
%A code named \texttt{UnitEvents} implemented for the \texttt{R} software 
%and accessible at \url{https://sourcesup.renater.fr/projects/uepackage} 
%has been developed to detect connectivity's between neurons from simultaneously recorded spike trains. The implemented procedure is based on a Gaussian approximation of the counts to determine coincidences between spikes. Some test procedures have also been developed to test if two neurons are independent or not.
%The present article is in line with these works.}

Regarding the membrane potential dynamic, it is commonly modeled through the Hodgking and Huxley model or FitzHugh-Nagumo neuronal model \citep[see \emph{e.g.}][]{samson}, using ionic channels data. When only the membrane potential is available over time, some simpler models are investigated, such as a diffusion process with jumps (\cite{JBHD}). This model focuses on catching the spiking phenomenon of the neuron itself. 

To the best of our knowledge, no paper studies the influence of spike trains from several neurons on the dynamics of the membrane potential of a fixed neuron, which is what we propose in this work.

\subsection{Main ideas and results}

In this work, we study the spike trains of $250$ neurons and the membrane potential of one central neuron. 
The last spike train of the $250$-size network is the one of the central neuron.
%We first extract the spikes from the membrane potential of the central neuron. 
From this large network of neurons, we want to get the interaction information. Indeed, the first goal of our study is to infer the connectivity graph and then focus on the interaction between the central neuron and the neurons around it. Neurobiologists usually get this information from empirical methods as cross-correlograms or PSTH (see \cite{pouzat}). We choose to investigate the spike trains using a multi-dimensional Hawkes process instead, to analyze the whole network simultaneously.  
The benefits of this modeling are plural: in particular, the multivariate aspect is central since it allows to correct for potential spurious correlations.  Moreover, we perform sparse inference procedures, which allow us to reconstruct a graphical model of neurons that interact together. 

Here we investigate two procedures designed for Hawkes processes to get the connectivity graph. The first one is called \texttt{ADM4} and is implemented in the Python Library \texttt{tick}. The algorithm is developed in \cite{ADM4}. \texttt{ADM4} is a parametric method based on the exponential assumption for the Hawkes process. It relies on a likelihood with a LASSO penalty criterion and is adapted to high-dimensional networks.
The second method is a nonparametric estimator proposed in \cite{lambert}.
This estimator is derived from the optimization of a penalized least-squared criterion and does not make any assumption on the interaction functions. However, this estimator was developed for a small network and to the best of our knowledge its performance has not been demonstrated in a high-dimensional context.

In Section \ref{sec:hawkes} we investigate the whole network with \texttt{ADM4} and then we extract a sub-network of $17$ neurons that are assumed to impact the central neuron (according to the outcome of the estimation of the adjacency matrix). We then re-estimate the parameters of both models from this subnetwork with both methods this time and compare the results.

Furthermore, we apply the test procedure introduced in \cite{reynaud2014goodness} on the subnetwork in order to validate the choice of the Hawkes process to model the interaction between neurons. This is an asymptotic test based on the Kolmogorov statistics, which is detailed in Section \ref{sec:test}.
%\textcolor{red}{à voir si on laisse ça ici :}
A code named \texttt{UnitEvents} implemented for the \texttt{R} software 
%and accessible at \url{https://sourcesup.renater.fr/projects/uepackage} 
has been developed to detect connectivities between neurons from simultaneously recorded spike trains. The implemented procedure is based on a Gaussian approximation of the counts to determine coincidences between spikes. Some test procedures have also been developed to test whether the spikes of two neurons are independent or not.
The present article is in line with these works.
We apply the test with the estimator developed in \cite{lambert} as in the original paper \cite{reynaud2014goodness}, but also with the estimator \texttt{ADM4}.
Let us note that this test is an asymptotic test with the number of samples and we only have ten samples at our disposal. However, since we observe many spikes on the study interval, we can assume that we are close to the asymptotic scenario and that the result remains interesting in our framework. 

In the last part of the article, we focus on a model for the membrane potential between spikes. We propose a jump-diffusion model with jumps driven by a Hawkes process. This model was studied theoretically in \cite{DLL}, \cite{DL2020} and \cite{ADGL2020}. The idea behind this model is that the membrane potential of a neuron can be modeled using a diffusion process, and we add the jumps driven by a Hawkes process to take into account the signals received by the central neuron from the neurons around. In this multivariate Hawkes process, we only keep the neurons interacting with the central neuron from the connectivity graph. We assume that these neurons impact the membrane potential of the central neuron. 
From the previous works \cite{DL2020} and \cite{ADGL2020} on this model, we know how to estimate all the coefficients non-parametrically. Then, we can generate new trajectories and compare them with the actual signal. 
 
Finally, we use the depth definition of Tukey \cite{tukey} adapted for curve data in \cite{pavlo} to see the importance of the self-exciting diffusion process to model the trajectories inter-spikes, instead of a simple diffusion process. 
 
The code (in \texttt{R} and \texttt{Python}) for the whole procedure is available on the repository \url{https://github.com/charlottedion/HawkesForNeuro}. Some details are given in Sections \ref{sec:resspikes} and \ref{sec:jumpdiffusion}.\\

%------
The paper is organized as follows. Section \ref{sec:data} is dedicated to the description of the data. We present both inference procedures for the study of the spike trains in Section \ref{sec:hawkes}. The inference results on the data, together with the result of the goodness-of-fit test, are given in Section \ref{sec:resspikes}.
In Section \ref{sec:jumpdiffusion} 
the jump-diffusion model is presented and implemented on the data. We present the estimation of the coefficients obtained on the data. With a simulation of new trajectories obtained from the inferred process, we obtain a validation of the model through the depth motion of Tukey.
Section \ref{sec:discussion} gathers the discussions on modeling and inference on this data.

\section{Data description}\label{sec:data}

\begin{figure}
\centering
\includegraphics[scale=0.4]{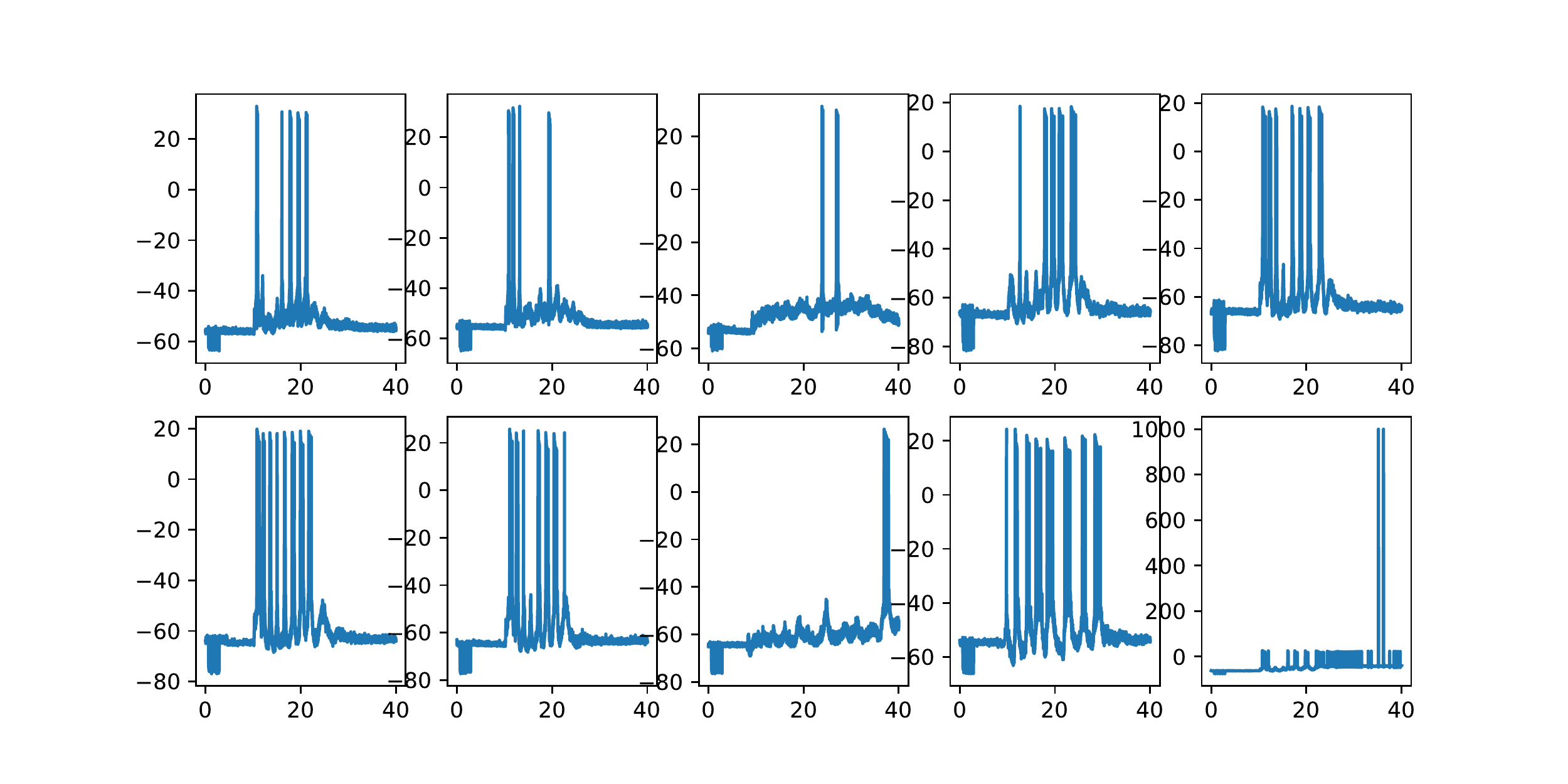}
\includegraphics[scale=0.4]{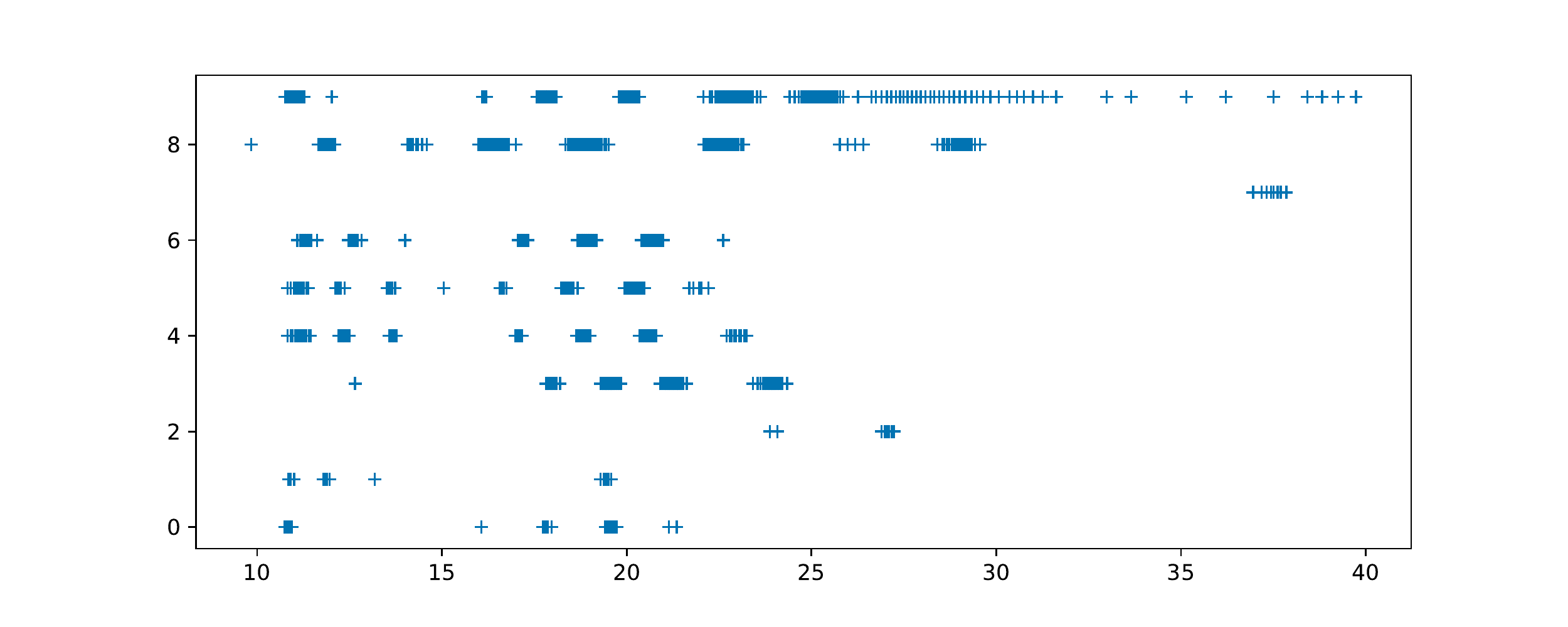}
%\caption{Extracted spike trains of the 10 trials of the fixed neuron}
\caption{Top graph: membrane potential of the 10 trials of the fixed neuron. Bottom graph: extracted spike trains of the 10 trials of the fixed neuron}
\label{fig:intraextra}
\end{figure}

The data consist of two kinds of signals from the lumbar spinal of a red-eared turtle: a continuous signal from the central neuron, the membrane potential of which has been recorded, and a discrete signal representing the spike trains of the neurons around. They are presented, for example, in \cite{radosevic}.
The two types of data were acquired with different recording systems and were synced using a trigger signal. They are measured simultaneously, and an external stimulus is applied to the turtle after $10$s during $10$s. 
%(signal shown in the intracellular recordings).
The first type of data is called
\textit{intracellular}.
It is an intracellular recording of the membrane potential of one neuron in mV during $40$s. The sample time step is $0.048$ms for the path. %and the observation interval is $[T^*, T_{max}]=[]$ thus there is $834000$ times for one path..\\
An invasive recording method allows getting this type of information that describes the observed neuron's electric dynamic. The experiment was conducted ten times. The ten paths are shown on Figure \ref{fig:intraextra}, upper graph. 

The second type of available data is called
\textit{extracellular}. The extracellular recording is also $40$ seconds long, and we observe $249$ neurons surrounding the central neuron. For each neuron, the recording is composed of the \textit{spike train} observed on the observation window. 
These $249$ discrete paths correspond to a network of neurons. Indeed, as they are observed in the same neuronal area, it is reasonable to think that they are communicating. 
This is discussed in detail in the following. 
Ten trials that we can consider as repetitions of the same signal are available.
On Figure \ref{fig:histallspikes}, we show on the left graph the spike trains for the first trial and on the right graph the distribution of the number of spikes for the first trials for the whole observations.

The first step in our analysis consists of recovering the connectivity graph between the neurons, especially the connectivity between the central neuron and the neurons around. We have to extract the spikes of the central neuron from its membrane potential recordings for that purpose.
We consider that the neuron spikes as soon as its signal exceeds $-20$mV. The results for the $10$ trials is presented of Figure \ref{fig:intraextra}, lower graph. 
%The spike trains of the central neuron for each of the 10 trials are represented on Figure \ref{fig:spikesX})
%Extraction of the spikes of central neuron. 
%Discussion on the threshold (short sensibility analysis).
We now have $250$ spike trains and $10$ trials. From Figures \ref{fig:intraextra} and \ref{fig:histallspikes}, we choose to focus on the time interval $[11,24]$ in the following (see Section \ref{sec:hawkes} for more details about this choice).
We remove $2$ neurons with no spike in the time interval of interest and the trial number $8$ for the same reason for the study presented below. We have in the following $M=248$ neurons and $n=9$ trials.

%\begin{center}
%\end{center}

\section{Inference of connectivity graph using Hawkes processes}\label{sec:hawkes}

In this section, we use temporal point processes to model the multiple spike trains. They are stochastic processes composed of time series of binary events that occur in continuous times \citep[see][]{DVJ}.
In particular, we consider Hawkes processes. We first present classical methods and explain why the automatic methods we use are more adapted and convenient. 

\subsection{Preliminaries}

Due to the external stimulation, the observed system is not stationary. Nevertheless, after the stimulation, it is possible to see a stability period. The Peristimulus Histogram (PSTH) is an existing tool to help with the decision about this phase.

\begin{figure}
\centering
\includegraphics[width= 7cm, height=5cm]{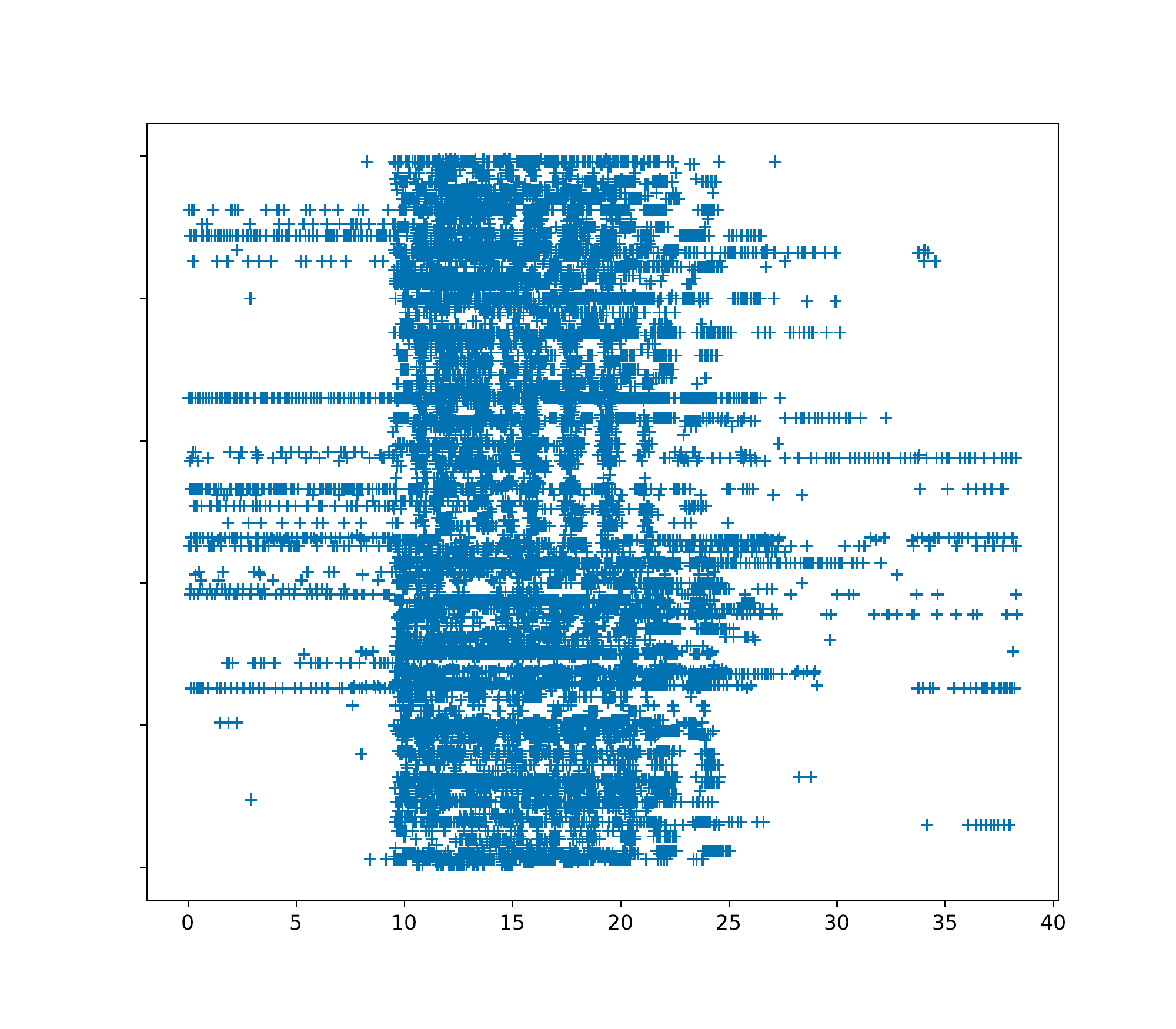}
\includegraphics[width= 7cm, height=5cm]{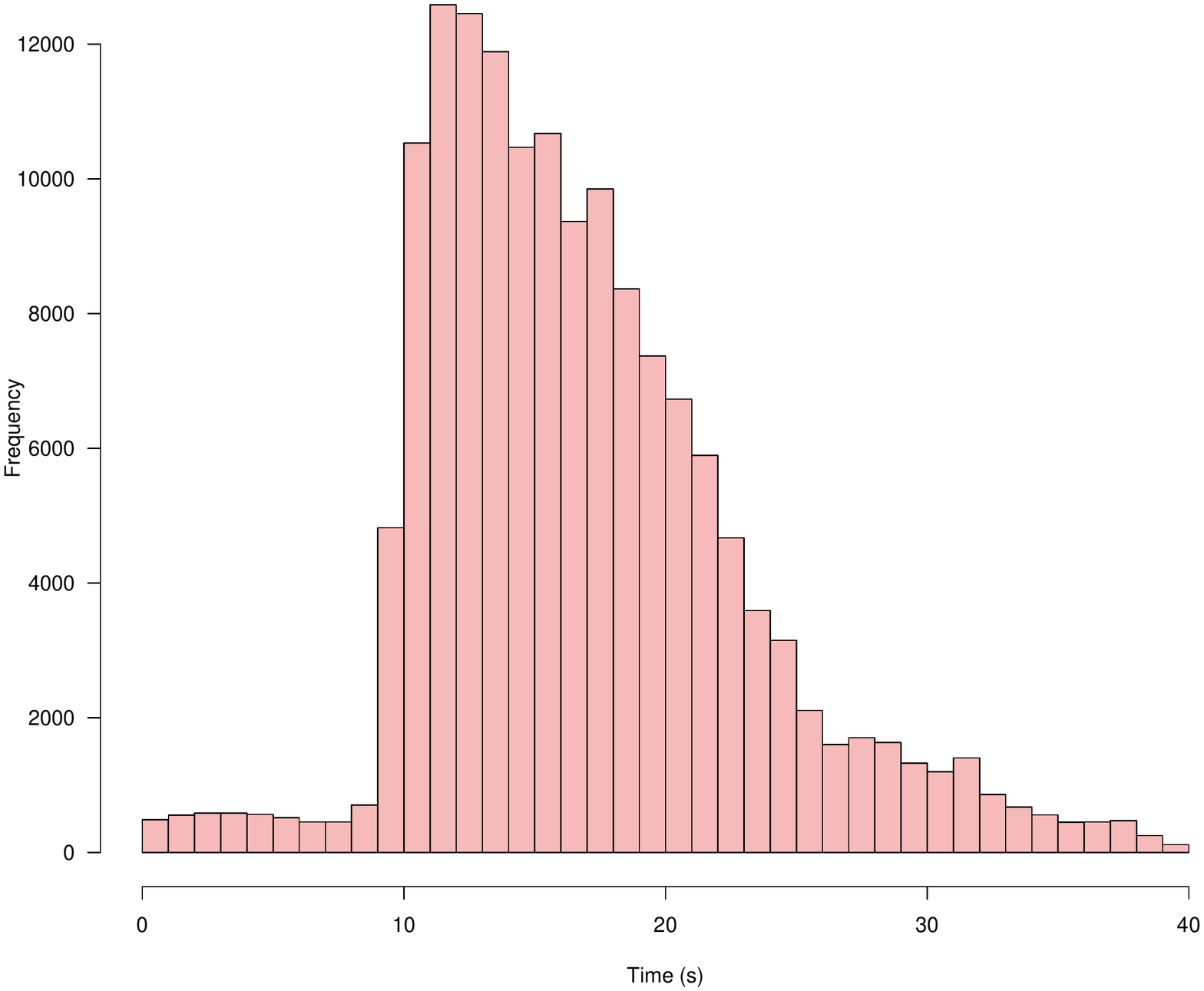}
\caption{Right: spikes trains for the first trials for all neurons. Left: distribution of the number of spikes for all neurons on the first trial (PSTH)}
\label{fig:histallspikes}
\end{figure}

\paragraph{PSTH for stationarity}
We have to find a time interval where we consider that the distribution of the points in a neighborhood of $t$ does not depend on $t$.
For this purpose, we consider the PSTH of all neurons to see if we can find a period during the stimulation phase where the signal can be assumed to be stationary. Indeed, the PSTH gives some measure of a neuron's firing rate or firing probability as a function of time, starting with (or around) a reference point (stimulus marker). In our data, the stimulation takes effect after $10$s of recording and lasts $10$s for each neuron and each trial.
We use the R-package \texttt{STAR}. The result is presented in Figure \ref{fig:histallspikes}. On the left, we have represented the spike trains for the first trial. On the right, the PSTH is presented. We choose to investigate the time interval $[11, 24]$s for the following inference.

\paragraph{Limits of classical tools}

Different methods are based on histograms of pairwise cross-correlation functions of neurons \citep{CCH}. From these procedures, some imprecise results may result in confounding neurons present in the network and not used in these histograms. The following methods of inference avoid this issue.

\subsection{Hawkes model methods}\label{sec:hawkesmethods}

%Before leading the study using the Hawkes model, one can investigate the Poisson process assumption. Much attention has been paid to  goodness-of-fit tests in this case in the technical report of \cite{pouzat} and in the theoretical work \cite{reynaud2014goodness}. 
%NONC'est nimporte quoi On naive visual test consist in represent the increments of counting process for each neuron. Indeed, in the Poisson probabilistic model, the consecutive increments are independent. 
% produces examples where the relation is clearly linear for neurons 114, 151, 249. \\
%-------

%with the selected neurons containing the central neuron. 
%The aim is to reconstruct the functional connectivity between the neurons. Then we focus on the neurons which impact the central neuron.
%Due to the large dimension of the dataset, the first goal is to reduce the number of neurons in the network. 
%Indeed, it is reasonable to think that the true number of neurons interacting are of order 10 \citep[][]{lambert,RBRTM}. 

%In this section we focus on the whole extracellular data.

We choose to model jointly the $M$ spike trains using a multivariate linear Hawkes process. This process is characterized by the conditional intensity of each coordinate (neuron) $\lambda_j(t)$:
\begin{equation}\label{eq:hawkes}
\lambda_j(t) = \mu_j + \sum_{j'=1}^M \sum_{k:T_k^{(j')}<t} h_{j,j'} (t-T_k^{(j')})
\end{equation}
where $\mu_j>0$ is the baseline intensity (or exogeneous intensity) and $h_{j, j'}$ is the interaction function (sometimes called synaptic weight function) of $j'$ on $j$, finally $(T_k^{(j')})_k$ is the sequence of jumps of neuron number $j'$ observed until time $t$.
% on $[T^*, T_{max}]$. 
The parameters to estimate are $\mu_j, h_{j,j'}, j=1, \ldots, M, j'=1, \ldots, M$.
Let us present the inference method used in the following. 

\paragraph{Parametric algorithm \texttt{ADM4}.}
We investigate the classical exponential kernel choice through the algorithm of \cite{ADM4} developed on the Python library \texttt{tick} under \texttt{HawkesADM4}.
%which proposes various methods appropriate for high dimension. 
%It is about the %\texttt{HawkesSumExpKern} and
The kernel functions are assumed to have the exponential decay shape as:
$$h_{j,\ell}(t)=a_{ j, \ell} \beta e^{-\beta t}$$
where $0<a_{ j, \ell}<1$ is the impact from neuron $\ell$ on neuron $j$.
The matrix of entries $(a_{j,\ell})$ is assumed to be sparse.
%\sout{ and low rank}.
 The estimation of the parameters is based on a penalized likelihood criterion with a LASSO penalty. The parameter $\beta$ is chosen as the best parameter in a grid in the sense of the least-squares criterion.

 This algorithm solves the estimation
problem efficiently by combining the alternating
direction method of multipliers and
minimization. It is a regularized convex optimization approach.

\paragraph{Nonparametric method with LASSO \texttt{NPL}.}
The second method we investigate is proposed in \cite{Hansen15} and \cite{lambert}. It is a nonparametric procedure, the main idea of which being to decompose the kernel functions on a piecewise-constant basis. More precisely, each kernel function $h_{j,\ell}(t)$ will be estimated by
$$g_{j,\ell}(t)= \sum_{k=1}^K \alpha_{j,\ell}^k \one_{((k-1)\delta, k\delta]}(t),$$

where the unknown parameters $\alpha_{j, \ell}^k$ are obtained by minimizing a least-squares contrast with LASSO penalty. 
 \cite{lambert} ensures that as long as the choice of the support $K\delta$ remains reasonable with respect to the excitatory or inhibitory behavior of the neurons, modifying the values of $K$ or $\delta$ does not change the reconstructed connectivity graph for a given value of the LASSO regularization parameter.

Both \texttt{ADM4} and \texttt{NPL} procedures have convergence results with the length of the time interval or with the number of trials (which have the same role using a concatenation of the trials). 
The convergence of the penalized least-squares with the same penalty used for \texttt{ADM4} is investigated in
\cite{BACRY}. In both cases, it is natural to think that the convergence can be achieved if the number of events is large enough. Besides, the observation time must be put into perspective according to the observed phenomenon; here, for example, $10$s is a large observation window for a neuron.

 \paragraph{Adjacency matrix}
The main interest is to reconstruct the following matrix $A$:
\begin{equation}\label{eq:A}
A= (a_{j,\ell})_{j,\ell\in \{1, \ldots, M\}}, ~ a_{j,\ell}= \int  h_{j,\ell}(t)dt
\end{equation}
% (can take negative values).
Each coefficient $a_{j, \ell}$ represents
 the average number of additional points created in the $j$th spike train  due to one spike on the $\ell$th spike train. It can be called the adjacency matrix for the graph where the $M$ neurons consist of the nodes. 
 In our case, the largest matrix is of size $M \times M$. 
 A non-zero coefficient in the matrix $A$ corresponds to an interaction between two neurons, making the adjacency matrix a convenient tool to visualize the connectivity graph between a group of neurons.

\paragraph{Comments}

Let us list the main advantages and drawbacks of both approaches. First, \texttt{ADM4} is 
easy to use (from library \texttt{tick} in Python) and adapted for large network size. The main drawback is that this algorithm is not developed yet for taking into account possible inhibition, meaning non-positive coefficient $a_{j, \ell}$. 
%\item[-] no inhibition
%\item[-] stability (res-sampling difficult) 

The second procedure \texttt{NPL}
has strong theoretical foundations and the nonparametric approach is less restrictive regarding the shape of interaction that can be non-exponential. Moreover, the estimation procedure can provide negative coefficients which would describe inhibitive effects. 
%Also, inhibition is allowed. 
The main weakness of the method is that it is adapted for small network size, as the authors have specified in \cite{lambert}.

We will study the relevance of using Hawkes processes to model the multiple spike train data. However, one can legitimately wonder the impact of the signals from unrecorded neurons, which is an important issue. One way to see it in the Hawkes model is the baseline parameters $\mu_j$, which are also called exogenous rates and can catch some of this hidden information.

%Indeed, in addition to what was said before, the dark side of the network is in a way taken into account in the model. It is an important issue to take into consideration the signals from unrecorded neurons. One way to see it in the Hawkes model is the baseline parameters $\mu_j$, which are also called exogenous rates and can catch this hidden information. 

\subsection{Validation of the Hawkes assumption}\label{sec:test}

Although the validation of the Hawkes assumption is crucial to ensure the relevance of the subnetwork we study, it is a very challenging aspect, particularly in the multivariate context. Indeed, a classical test for point process model is based on the so-called \textit{time-rescaling theorem}, presented for example in \cite{brown2002time, embrechts} and proven in \cite{papa}, but it raises some issues highlighted in  \cite{reynaud2014goodness}. These issues come from the fact that we use an estimator of the intensity in the test, which has been obtained from the same data used in the test. To overcome these issues, \cite{reynaud2014goodness} develops a test using sub-sampling together with cumulative use of the trials.

Let us present this test (it corresponds to Test 4 in \cite{reynaud2014goodness}). 
For the $M$ neurons, we define by $N_j^{(i)}, i=1,\dots,n$ and $j=1, \ldots, N$ the spike trains of neuron $j$ for trial $i$ during the time interval $[0,T_{\max}]$. 
The $n$ trials are assumed \textit{i.i.d.}

The first step is to compute an estimator
 $\w{\lambda}^{(i)}_j, j=1, \ldots, M$ of the intensity for each sample $i$ (using \texttt{ADM4} or \texttt{NPL}).
 %The plug-in step leads to bad results as explained in \cite{reynaud2014goodness}.
Then, for each neuron $j \in \{1, \ldots, M\}$, execute the following steps: 
\begin{enumerate}
\item Take a sub-sample of trials denoted $\mathcal{S}$ (in $\{1, \ldots, n\}$) of size $p_n= [n^{2/3}]$.
\item For each trial $i$ in $\mathcal{S}$ change time and compute $\w{\Lambda}^{(i)}_{j}(T_k)= \int_0^{T_k} \w{\lambda}_j^{(i)}(s)ds$ for 
$T_k$ the jumps of the neuron $j$ in the trial $i$.
\item Cumulate\footnote{see \cite{reynaud2014goodness} for the definition of cumulative process. It basically means that the sequences are put together.} the $\w{\Lambda}^{(i)}_j$ processes for $i \in \mathcal{S}$, denote $\w{N}_j^{c,\mathcal{S}}$ the resulting process on $[0, \sum_{i \in \mathcal{S}}\w{\Lambda}^{(i)}_j(T_{\max})]$.
\item Fix $\theta >0$, strictly smaller than $\sum_{i \in \mathcal{S}}\w{\Lambda}^{(i)}_j(T_{\max})/p_n$.
\item Compute 
$$Z_j = \sqrt{\w{N}_j^{c, \mathcal{S}}([0, p_n\theta])} \sup_{0 \leq u \leq 1} \left| 
\frac{1}{\w{N}_j^{c, \mathcal{S}}([0, p_n\theta])} \sum_{X \in \w{N}_j^{c, \mathcal{S}}, X \leq p_n\theta} \one_{X/(p_n\theta) \leq u }-u
 \right|$$
 \item Reject if $Z_j > q_{1- \alpha}$.
\end{enumerate}
Here $q_{1- \alpha}$ is the asymptotic quantile of the Kolmogorov-Smirnov distribution. In the tables, we get the empirical quantile $q_{p_n,1-\alpha}$, we deduce the asymptotic quantile by considering $(\sqrt{p_n}+0.12+0.11/\sqrt p_n)\times q_{p_n,1-\alpha}$
%(see \cite{reynaud2014goodness}) or $(\sqrt{n}+0.12+0.11/\sqrt n)*q_{n,1-\alpha}$ 
as $p_n=4<45$ here (see \cite{stephens1974edf} for details).
%Remarks:
%\begin{itemize}
%\item[-] Asymptotic result $p=p_n \rightarrow \infty$, we have $p= 4$
%\item[-] Table of KS distribution.
%\end{itemize}
We drive the reader attention on the fact that this is an asymptotic result on $p_n \rightarrow \infty$ when $n \rightarrow \infty$, but in the present context we apply the test with $p_n= 4$ and $n=9$. Then we consider  the outcome of this test cautiously due to the small sample size.

However, let us precise that since we observe lots of spikes in the considered time interval, we can assume that we are not so far from the asymptotic framework and it seems reasonable to consider this test in our case.

In the next Section, we present the results regarding the Hawkes inference on the neuronal data.

%%%%%%%%%%%%%%%%%%%%%%%%%%%%%%%
\section{Study of the spike trains}\label{sec:resspikes}

This section aims to use the extracellular data to reconstruct a connectivity network containing neurons that interact with the central neuron.

We first estimate the adjacency matrix $A$ on the whole network of $M=248$ neurons from the $n=9$ trials using \texttt{ADM4}. We extract from this estimation procedure a sub-network of $M=17$ neurons that impacts the central neuron. Then we perform the goodness-of-fit test with the estimated coefficients obtained from the two inference procedures: \texttt{ADM4} and \texttt{NPL}.

\subsection{Estimation of a large adjacency matrix}

As explained in Section \ref{sec:data}, we work with the spikes of $M=248$ neurons. Due to the large dimension of the matrix $A$ given in \eqref{eq:A}, we use the \texttt{ADM4} method to estimate it.
The decay parameter $\beta$ is estimated at $20$Hz.
Only $46.51\%$ of the coefficients are non-zero. 

We then focus on the last line on the estimated matrix $\w{a}_{M,\ell}, \ell=1, \ldots, M$ which represents the impact coefficient of each neuron on the central neuron (number $M$).
We keep the neurons for which the coefficient is strictly positive (larger than $10^{-5}$) which, according to the model, means that the dynamic of the intensity of the central neuron is increased when one of these neurons spikes. We find $17$ neurons given in Table \ref{tab:tablecompar}. In column "\texttt{ADM4}" we give the corresponding estimated coefficient $\w{a}_{M, \ell}$ in the large network. 

\begin{table}
    \centering
    \begin{tabular}{c|cc}
     \hline
    \hline
    Neuron &  \texttt{ADM4} & Empirical \\
    \hline 
 48  & 0.003 & 0.010 \\
 58  & 0.027 & 0.056\\
 72  & 0.085 & 0.031\\
 99  & 0.088 & 0.012\\
 115 & 0.006 & 0.006\\
 116 & 0.039 & 0\\
 119 & 0.015 & 0.020\\
 126 & 0.043 & 0.034\\
 140 & 0.028 & 0.018\\
 152 & 0.03  & 0.020\\
 153 & 0.081 & 0.067\\
 192 & 0.04  & 0.017\\
 194 & 0.068 & 0 \\
 210 & 0.079 & 0.061\\
 230 & 0.165 & 0.146\\
 232 & 0.025 & 0.029\\
 240 & 0.059 & 0.030\\
  \hline
    \hline
 \end{tabular}
 \caption{\texttt{ADM4} estimation obtained on the large network and empirical coefficient obtained from the 9 trials for the $17$ neurons of interest}
 \label{tab:tablecompar}
\end{table}

\paragraph{Comparison with spike-triggered membrane potential method}

We take advantage of the available intracellular recording to proceed to a backward validation. This
empirical method, suggested by \cite{petersen}, allows detecting spike-triggered membrane potential. 
The idea is to validate of the assumption that the neurons from the estimated subnetwork and the central neuron interact. To this purpose, we investigate whether we can observe a change in the membrane potential of the central neuron when the other neurons spike. More precisely, we compute for each neuron $j$ a coefficient which is the ratio between the number of spikes of the central neuron on each short interval around each spike of the neuron $j$ (from 2ms before until to 2 ms after).  We display in Table  \ref{tab:tablecompar} column "Empirical" this coefficient averaged on all spikes of each neuron $j$. We can observe that the largest values of the empirical coefficients correspond to the positive values of the \texttt{ADM4}. %coefficients which comforts us .
In the following, we pursue the study with this selected subset of neurons.
%based on these two estimations.

% For each spike of the neurons of the sub-network that are supposed to impact the membrane potential of the central neuron, the corresponding phase of the continuous signal. 
%On these phases of $4$ms ($2$ms before and $2$ms after spike as it is done in \cite{petersen}) we are looking for an impulsion. The empirical method is illustrated on Figure~\ref{fig:plotrune}.
%Then we compute for each neuron $j$ a coefficient which is the ratio between the number of spikes of the central neuron visible on the membrane potential on each short interval around a spike of the neuron $j$.

%More precisely, we display these coefficients on Table \ref{tab:tablecompar} first column. We can observe that the largest values of the Empiric coefficients correspond to the positive values of the \texttt{ADM4} coefficients which is really comforting. We can continue now studying this small network based on these two estimations.

%\begin{figure}
%    \centering
%    \includegraphics[width=10cm, height=7cm]{plotrune.pdf}
%    \caption{Illustration of spike triggered method, for neuron 230 on trial 4}
%    \label{fig:plotrune}
%\end{figure}
%\textcolor{blue}{Figure \ref{fig:plotrune} non commentée}

\subsection{Inference on the sub-network}
\label{sec:subnetwork}

%Line one is the result on the time interval $[11, 24]$, second line on $[11, 20]$.

%{\small \begin{table}[H]
%    \centering
 %   \begin{tabular}{cccccccccccccccccc}
%48 & 58&  72&  99& 115& 116& 119& 126& 140& 152& 153& 192& 194  &210& 230& 232& 240\\
%48 & 58&  72&  99& 115& 116& 119& 126& 140& 152& 153&& 192& 194&  &210& 230& 232& 240\\
%& & 72 & 99 &106 &118& 119& 126& 140 & && 162& 192&  &206 &210& 230& 232 &
%    \end{tabular}
 %   \caption{Selected neurons impacting the central neuron. On the time interval $[11, 24]$.}
  %  \label{tab:selectneurons}
%\end{table}
%}

\begin{figure}
    \centering
    \includegraphics[scale=0.8]{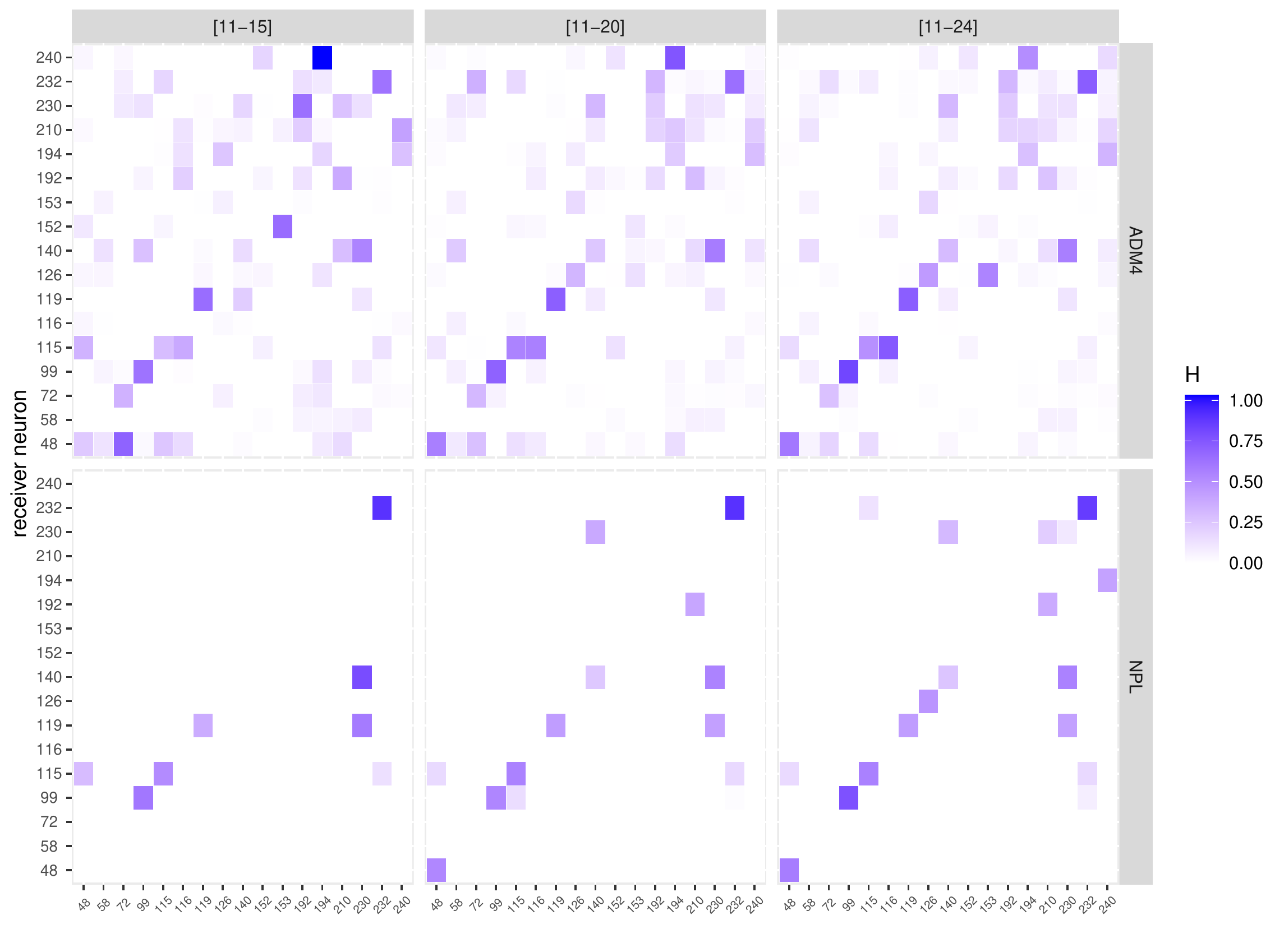}
    \caption{Estimated adjacency matrix $A$ on the network of size $17$ with \texttt{NPL} (bottom) \texttt{ADM4} (up) on different time intervals. }
    \label{fig:comparison_ADM4_NPL_Tmax}
\end{figure}

On the sub-network of size $M=17$, which extraction procedure has been described before,
we can use both estimation procedures described in Section \ref{sec:hawkesmethods}: \texttt{ADM4} and \texttt{NPL}. The result for the estimation of $A$ is given in Figure \ref{fig:comparison_ADM4_NPL_Tmax}: bottom line for \texttt{NPL} and top line \texttt{ADM4}. The estimation is represented as a function of the time observation interval. There are three graphs for three different periods of observation considered: $[11, 15], [11, 19], [11, 24]$. In the following study, we consider the results obtained on the largest interval since it contains more occurrences. However, the interest of this dynamic representation is to investigate the stability in time of the estimations with both inference procedures. This stability is visually verified in Figure \ref{fig:comparison_ADM4_NPL_Tmax} since there are very small variations according to time in the estimated adjacency matrices. In table \ref{tab:compar_ADM4_NPL_Tmax_bis}, we have computed the Frobenius norm and the spectral norm of the difference between the adjacency matrices obtained from \texttt{ADM4} and \texttt{NLP}. Interestingly, Table \ref{tab:compar_ADM4_NPL_Tmax_bis} shows that the larger the observation interval (that is, the more spikes observed), the more similar the adjacency matrices. We also notice that \texttt{NPL} provides estimated matrices which are sparser than those obtained with \texttt{ADM4} (see Table \ref{tab:compar_ADM4_NPL_Tmax}). Nevertheless, the largest coefficients have similar estimated values with both methods, suggesting that the two inference procedures provide consistent results. 
%\textcolor{red}{Changer les commentaires en dessous + commenter les deux tables que l'on garder et vérifier les résultats}
%We can see that the estimation is more sparse with \texttt{NPL}. Nevertheless the main coefficients have the same estimated values.
%The matrix obtained with \text{ADM4} has $61.5\%$ of non-zero coefficient versus $17.3\%$ for \texttt{NPL}.

%Some comparison elements are given in Table \ref{tab:estimAM17} where \texttt{Fro} for the Frobenius norm and 
%$\rho$ is the spectral radius. 

%Besides, let note that the estimation obtained with \texttt{ADM4} seems robust in the sense that the estimation obtained on $[11, 15]$ and $[11, 20]$ are very close to the one obtain on the entire interval $[11, 24]$. Nevertheless the estimation seems to be more precise with more points.

\begin{table}
    \centering
    \begin{tabular}{c|ccc}
     \hline
    \hline
    Estimator & [11,15]&  [11, 19] & [11,24]\\
   \hline
         \texttt{ADM4}& 56.4 & 58.8 & 57.4\\
        \texttt{NPL}& 2.8  & 4.8 & 6.2\\
         \hline
    \hline
    \end{tabular}   
    \caption{Percentage of non-zero coefficients in the estimated adjacency matrices}
   \label{tab:compar_ADM4_NPL_Tmax}
\end{table}

\begin{table}
    \centering
    \begin{tabular}{c|ccc}
     \hline
    \hline
    Norm & [11,15]&  [11, 19] & [11,24]\\
   \hline
         Frobenius norm & 2.14 & 1.59 & 1.46\\
      Spectral norm & 1.13  & 0.91 & 0.76\\
       \hline
    \hline
    \end{tabular}   
    \caption{Norms of the difference between adjacency matrices estimated with ADM4 and NPL}
   \label{tab:compar_ADM4_NPL_Tmax_bis}
\end{table}

%\begin{figure}
 %   \centering
  %  \includegraphics[scale=0.15]{adm-pat-17.jpeg}
   % \caption{Estimated adjacency matrix $A$ on the network on size $17$: left \texttt{NPL}, right \texttt{ADM4}. }
    %\label{fig:comparisonA}
%\end{figure}

\paragraph{Goodness of fit test}

We apply the procedure described in Section \ref{sec:test} to test the assumption that the multivariate process composed of the spike trains of the $17$ neurons previously selected is a Hawkes process. The null hypothesis here is that the data comes from a Hawkes model, and we perform the test with a standard risk level $\alpha= 0.05$. Let us mention that the method is described for each neuron, and for one subsample of randomly chosen trials.
For the neurons $119, 153, 192$ for example, we can notice that the estimations close (see Figure \ref{fig:comparison_ADM4_NPL_Tmax}) and the obtained rates are also close.
Regarding this last point, we propose here to perform the procedure for $100$ randomly chosen subsamples to increase its robustness. We compute an empirical acceptance rate for each neuron, which corresponds to the percentage of times the null hypothesis has been accepted over the $100$ tests. We provide the outcomes of this goodness-of-fit test in Table \ref{tab:Ztesttable}. Although some acceptance rates are pretty low, most of them are greater than $50\%$ for both \texttt{ADM4} and \texttt{NPL} procedures. As a comparison, we generated a dataset according to an exponential Hawkes process with an adjacency matrix that corresponds to the one estimated with \texttt{ADM4}: the results are denoted as "synthetic data" in Table \ref{tab:Ztesttable}. We observe that, as expected, the acceptance rates are high on average, but even on this synthetic dataset, some acceptance rates are lower than $50\%$, which comforts us in validating the Hawkes assumption on our neuronal dataset.

\begin{table}[ht]
    \centering
    \begin{tabular}{c|ccc}
    \hline
    \hline
    Neuron & \texttt{ADM4} &  \texttt{NPL} & Synthetic data  \\
    \hline
        48 & 0.21 & 0.72 & 0.13\\
        58 & 0.46 & 1 &0.88  \\
        72 & 0.79  & 1&0.86 \\
         99& 0.64 & 1&0.55 \\
         115 & 0.65 & 0.32 & 0.48\\
           116& 0.64 &1 & 0.89 \\
            119& 0.38 &0.44& 0.31\\
             126& 0.44 &0.1& 0.66\\
              140& 0.12 & 0.14&0.36\\
              152 & 0.27 &1& 0.7\\
               153 & 0.98 &1& 0.86 \\
                192 & 0.78 &0.81& 0.88\\
                  194& 0.89 &0.61& 0.43\\
                  210 & 0.28 &0.67& 0.5\\
                   230 & 0.43& 0.44&0.81  \\
                    232 & 0.04 & 0.51&0.39\\
                     240 &0.59  & 0.35&0.37\\
                      \hline
    \hline
    \end{tabular}
    \caption{Empirical acceptation rate for the goodness-of-fit test under the multivariate Hawkes process assumption for the real data using \texttt{ADM4} and \texttt{NPL}, and for the "synthetic data" simulated from ADM4}
    \label{tab:Ztesttable}
\end{table}
%%%%% Remarque: table pour theta = la somme *0.5

%{\small
%\begin{table}[ht]
 %   \centering
  %  \begin{tabular}{c|ccccccccccccccccc}
    %Neuron & 48 & 58 & 72 & 99 & 115 & 116 & 119 & 126 & 140 & 152 & 153 & 192 & 194 & 210 & 230 & 232 & 240 \\
   %  \texttt{NPL}& 0.71 & 0.36 & 0.69 & 0.71  & 0.94 & 0.76 & 0.79 & 0.61 & 0.63 & 0.88 & 1.0 & 0.65 & 0.76 & 0.57 & 0.79 & 0.04 & 0.93 \\
    % \texttt{ADM4} & 0.90 & 0.07 & 0.22 & 0.60 & 0.96 & 0.03 & 0.77 & 0.90 & 0.90 & 0.25 & 0.96 & 0.92 & 0.87 & 0.94 & 0.86 & 0.94 & 0.85  
     %    \end{tabular}
    %\caption{Empirical acceptation rate for the goodness-of-fit test under the multivariate Hawkes process assumption}
    %\label{tab:Ztesttable}
%\end{table}
%}

%%%%%%%%%%%%%%%%%%%%%%%%%%%%%%%%%%%%%%%%%%%%%%%%%%%%%
\section{Diffusion process to study the membrane potential}\label{sec:jumpdiffusion}

In this section, we consider the intracellular signal, which is a membrane potential recording. The aim is to find a model for the dynamic of this path between the spikes, taking into account the extracellular signal coming from the neurons around. 

\subsection{Self-exciting jump-diffusion model}

We rely here on the methodology developed in \cite{DL2020} and \cite{ADGL2020}. Let us assume that the following equation describes the membrane potential dynamics along time between spikes:
\begin{eqnarray}\label{eq:jumpdiff}
dX_t&=& b(X_t)dt + \sigma(X_t) dW_t+ a(X_{t-})\sum_{j=1}^M dN_j(t) ,~X_0=x_0.
\end{eqnarray}
%to describe the inter-spike-interval of a membrane potential trajectory.
%time interval: $[11, 16]$
The coefficient $b$ is the drift function, $\sigma$ is called the diffusion coefficient, and $a$ is the jump coefficient function. Note that $W$ is a standard Brownian motion independent of $N$, which is a $M$-dimensional linear exponential Hawkes process. The intensity $\lambda$ (given in (\ref{eq:hawkes}) with the exponential kernel) of the Hawkes process do not depend on $X$.

This process aims at describing how the dynamic of the membrane between spikes is impacted by the activity of the network of neurons surrounding one fixed neuron.
Then here, $M=17$ is the number of neurons that were detected impacting the central neuron (the neuron for which the potential is recorded) in Section \ref{sec:subnetwork}. From \cite{ADGL2020}, we know how to estimate $\sigma^2$ and $a$ thanks to the wealth of the present dataset. Then, from \cite{DL2020}, we can deduce an estimator of the drift. For the sake of simplicity, the estimators are not given here. We recall the definitions of the different estimators in the Appendix, Section \ref{sec:appendix}. We refer the reader to the original papers for more details concerning the definition of the estimators and the theoretical results. The code with the implementation of the three estimators $b$, $\sigma^2$ and $a$ is available on the GitHub page.

From this model, with the estimations of the parameters $b$, $\sigma^2$ and $a$, we are able to recover the trajectories of the membrane potential of the central neuron between spiking phases. In the following, we describe the procedure more precisely. 

\subsection{Description of the algorithm}

The inference is decomposed into several steps. We have $M=17$ neurons that interact with the central neuron. 
First, the estimations of $a$ and $b$ depend on the estimation of the intensity $\lambda$ of the Hawkes parameters. Thus we estimate $\lambda$ as described above with \texttt{ADM4}, and let us denote $\w{\lambda}$ the estimator. Then, we estimate $\sigma^2$ and $a$ from the procedures detailed in \cite{ADGL2020}. Finally, we recover the drift function $b$ from the method developed in \cite{DL2020}. More precisely, let us detail the main steps in the following.

\begin{enumerate}
    \item Data used for the inference on Model \eqref{eq:jumpdiff} (see Figure \ref{fig:datajumpdiff}): 
    \begin{itemize}
        \item focus on trial 1 for both the extra and intracellular signals,
        \item from the extracellular signal, keep the spike trains $N_j([T_*, T^*]), j=1, \ldots, M$ on the interval $[T_*, T^*]=[11, 16]$ where the process seems stationary,
        \item re-scale the spike trains on $[0,5]$,
         \item from the intracellular signal, we keep the membrane potential from trial 1 in the same time interval.
        \end{itemize}
    %\item Estimation of the intensity $\lambda_j$ for $j\in\{1,\dots,M\}$: from the spike trains of the sub network of size $M$ of the neurons impacting the central neuron, estimate $\lambda_j, j=1, \ldots, M$ from the ADM4 procedure: we get $\w{\lambda}_j$.
    
    \item Inference on $\sigma^2$ and $a$ following the procedure given in \cite{ADGL2020}: we obtain $\w{\sigma}^2, \w{a}$, then we approximate the obtained estimation by a constant for $\sigma^2$ and by a linear function for $a$.
    \item Inference on $b$ from the method develop in \cite{DL2020}: depending on the estimations $\hat\sigma^2$, $\hat a$ and also  $\w{\lambda}_j$, for $j=1,\dots, M$, estimate $b$ and obtain a realization of $\w{b}$, we also apply here a linear approximation.
\end{enumerate}

Then, from these estimators, we are able to simulate new trajectories that we can compare to the real trajectory of the membrane potential of $X$ on $[0,13]$ (we have rescaled the time interval $[11,24]$ so that it begins at 0). Let us notice that the time interval is longer than the one used in step 1 to find the connectivity graph. Indeed, here we consider a phase without spikes for the fixed neuron, and this phase is longer than the one used for the stationary of the extracellular signal. 
Let us detail the procedure: 
\begin{enumerate}
    \item From $\w{\lambda}_j$, $j\in\{1,\ldots,M\}$, simulate new spike trains $\tilde{N}_1, \ldots, \tilde{N}_M$ on $[0, T]$ with $T=13$. 
    \item Define functions $\w{b}, \w{\sigma}, \w{a}$ (which are linear approximations of the obtained estimators).
    \item Simulate trajectories $(X_{k\Delta})_{k\Delta \in [0, T]}$ from equation 
   \begin{equation}\label{eq:jumpdiffestim}
   dX_t= \w{b}(X_t)dt + \w{\sigma}(X_t)dW_t+ \w{a}(X_{t-})\sum_{j=1}^M d\tilde{N}_j(t)
   \end{equation}
   where the $\tilde{N}_j$'s are given by step 1. 
\end{enumerate}

\subsection{Results}

\begin{figure}
    \centering
    \includegraphics[width=10cm, height=7cm]{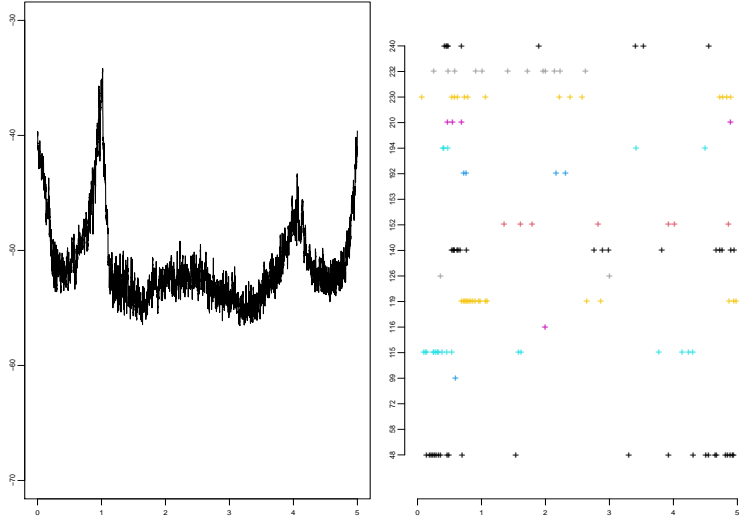}
    \caption{Representation of the data, left: the membrane potential of the central neuron, right: the spike trains of the 17 selected neurons that interact with the central neuron}
    \label{fig:datajumpdiff}
\end{figure}

The data used for this part are represented on Figure \ref{fig:datajumpdiff}: on the left, the section of the membrane potential that we want to describe, rescaled on $[0,5]$ and on the right, the corresponding spike trains for the $M=17$ neurons on the same interval. Let us mention that we have only kept one observation over ten for the intracellular signal because the data are too large otherwise. 

\paragraph{Coefficients estimation}

\begin{figure}
    \centering
    \includegraphics[width =14cm, height=6cm]{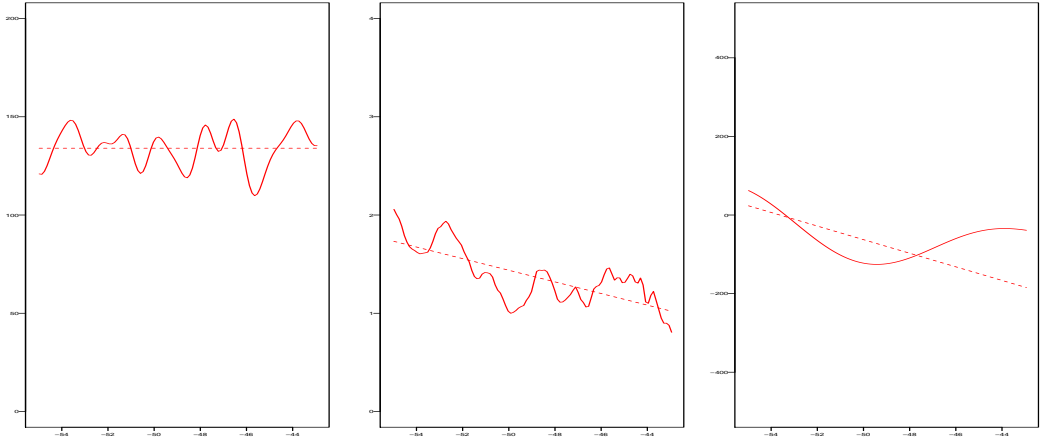}
    %{graphestim-ggplot.pdf}
%    {estim3coeff.pdf}
    \caption{Estimators of the parameters of the diffusion model with jumps, left: $\hat\sigma^2$, middle: $\hat a$, right: $b$; the plain red line represents the estimator and the dotted red line the approximation we made (constant for $\hat\sigma^2$, linear for $a$ and $b$)}
    \label{fig:estimcoeff}
\end{figure}

On Figure \ref{fig:estimcoeff}, we show in plain red line the estimated coefficients: $\w{\sigma}^2$, $\w{a}$ and $b$. We add on each graph the linear approximations we made in dotted red. Indeed, on the real data the estimators are really noisy and thus, we make approximations. Thus the coefficient $\sigma$ is assumed constant and $a,b$ linear. 

\begin{figure}
    \centering
    \includegraphics[width =14cm, height=6cm]{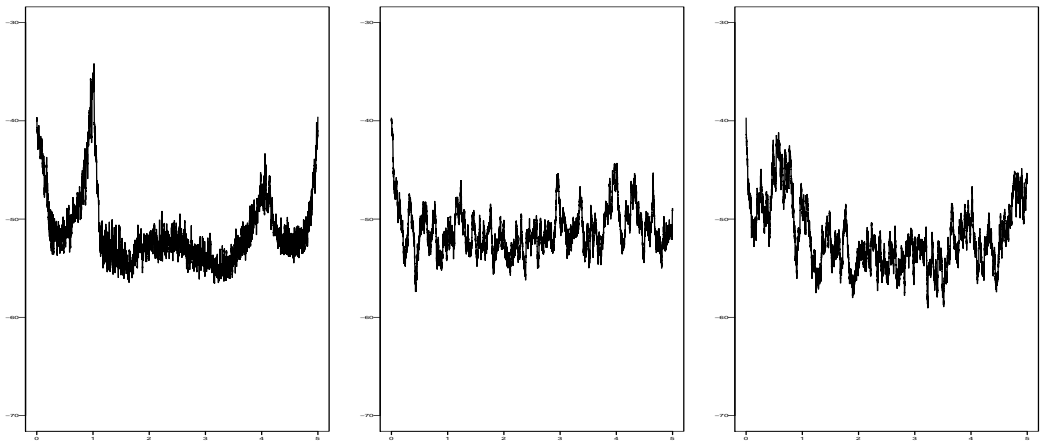}
    %{example-avecssaut.pdf}
    \caption{Trajectories of the membrane potential. Left: the true membrane potential (we have taken one observation over ten to represent the trajectory), middle: from a pure diffusion process without jumps, right: from model (\ref{eq:diffusion})}
    \label{fig:trajsimu}
\end{figure}

\paragraph{Inference path}
In Figure \ref{fig:trajsimu} we compare three trajectories:  the true portion of the membrane potential (right), the trajectory of $X$ simulated from model \ref{eq:jumpdiffestim} (right) and the trajectory of $X$ if we fit a pure diffusion process from the intracellular data without considering the spikes from the other neurons assuming there are no jumps (middle graph).

%represents on the third graph the true portion of the membrane potential, on the middle graph, we represent a trajectory of $X$ simulated from model \ref{eq:jumpdiffestim}. 
%Finally, on the first graph, we represent the trajectory of $X$ if we fit a pure diffusion process from the intracellular data without considering the spikes from the other neurons (assuming there are no jumps).
Precisely, it is possible to use a simple homogeneous process of type $dX_t= {b}(X_t)dt+ {\sigma}(X_t)dW_t$, to infer the trajectory without taking into account the information of the neurons around. 
In this case $b, \sigma^2$ are nonparametrically estimated with the procedure described in \cite{CGCR}, we denote $\wt{b}$ and $\wt{\sigma}$ the estimators.
We simulate then a trajectory through the model 
\begin{equation}\label{eq:diffusion}
dX_t= \wt{b}(X_t)dt+ \wt{\sigma}(X_t)dW_t.
\end{equation}
The result obtained in Figure \ref{fig:trajsimu} seems visually convincing that considering the spike trains of other neurons can improve the prediction of the membrane potential trajectory. Indeed, the middle trajectory simulated from the process with jumps seems closer to the true trajectory than the one simulated without spikes. In the next paragraph, based on depth notion, we propose a quantitative argument to support this observation.

\paragraph{Depth for curve}

In order to assess the interest of accounting for jumps in the diffusion model, we propose to use the notion of depth of Tukey developed for functional data in the R-Package \texttt{curveDepth} of  Pavlo Mozharovskyi, based on the results of \cite{pavlo}. The idea is to evaluate the similarity between paths. When the depth coefficient (which belongs to $[0,1]$) is large, the trajectories have many similarities. On the contrary, when it is close to zero, the paths are likely not to come from the same model.
The function \texttt{depthc.Tukey} allows us to compute the depths for the two simulated paths (with and without jumps) and compare them with the true ones.
Let us detail the procedure:
\begin{enumerate}
    \item Simulate a sample $D^{H}_{N_{\rm rep}}=(\w{X}_t^i)_{t\in[0,T]},i=1, \ldots, N_{\rm rep}$ from Equation (\ref{eq:jumpdiffestim}) 
    with $X_0 \sim \mathcal{U}({[-55, -35]})$. Simulate a sample $D^{wtH}_{N_{\rm rep}}=(\wt{X}_t^i)_{t\in[0,T]},i=1, \ldots, N_{\rm rep}$ from Equation (\ref{eq:diffusion}), with $X_0 \sim \mathcal{U}({[-55, -35]})$.
    \item Compute the depth of each curve of the sample $D^{H}_{N_{\rm rep}}$ (respectively $D^{wH}_{N_{\rm rep}}$) on the sample $D^{H}_{N_{\rm rep}}$ (resp $D^{wH}_{N_{\rm rep}}$).
    \item Compute the depth of the real trajectory $(X_t)_{[0,T]}$ on each sample and save the ranking. 
\end{enumerate}
We choose $N_{\rm rep}= 100$ and perform $50$ Monte-Carlo repetitions for this procedure. 
The depth of the real trajectory on the sample simulated from process (\ref{eq:jumpdiff}) is $0.67$ and the ranking is $60$ over $100$. The high rank of the true observation on the sample $D^{H}_{N_{\rm rep}}$ suggests many similarities between the generated and the true trajectories.
The depth of the real trajectory on the sample simulated from a simple diffusion process is $0.44$ and the ranking is $1$ over $100$. 
This study validates using a Hawkes-jump-diffusion process on this data and the closeness of the true data to the path simulated under the model \ref{eq:jumpdiff}.

\section{Conclusion and discussion}
\label{sec:discussion}

In this paper, we propose a complete study of real neuronal data. We integrate both intracellular and extracellular signals jointly in a unified framework to represent the dynamic of the membrane potential of a neuron considering the spikes from other neurons around. We first extracted a subnetwork of neurons that interact with the neuron of interest. We then studied the relevance of considering a Hawkes process to model the interactions between neurons based on a validation test. Finally, Model (\ref{eq:jumpdiff}) was introduced to take into account both types of data. Theoretical results to guarantee the performance of the proposed estimators were proposed in \cite{DLL}, \cite{DL2020} and \cite{ADGL2020}. In this paper, we work with this model on neuronal data to investigate whether it could improve the fit to the real data compared to a standard diffusion process. We compare both models (with and without jumps) using the depth of the curves representing the trajectories and we show that Model (\ref{eq:jumpdiff}) improves the fit to the real data.  In conclusion, our work highlights that even if we can only measure the spike trains of a limited number of neurons around the central neuron, taking this information into account in the model improves the fit of the membrane potential dynamics between spikes.

This work is a first step in the joint study of extra and intracellular signals from neurons. The Hawkes-diffusion can be further improved by considering, for example the jumps of the neuron itself or by taking a different a-function for each neuron that interacts with the central neuron. This work is beyond the scope of the paper, but it would be interesting to see if these perspectives would allow us to improve the fit of the trajectories. 

%Model (\ref{eq:diffusion}) that we propose to study 
Additionally, let us note that the inhibition behavior of the neurons
could be taken into account in the Hawkes-diffusion process. Indeed, the theoretical paper \cite{DLL} shows the results for a non-linear kernel with possible inhibition. This would require an efficient method to estimate the interaction functions of an exponential Hawkes model in a context of potential inhibition. A maximum likelihood estimator has been recently proposed in \cite{BMS21} for a univariate process and the extension to the multivariate case is an ongoing work. 
Finally, the last issue is to propose an efficient method that can handle the large dimension of the network.
For example the recent work \cite{newref} investigates large networks adapting  the method \texttt{NPL}. This is a serious lead for further works. 
%Nevertheless, for the diffusion part the exponential kernel is crucial to keep markovianity of the coupled process.

\paragraph{Acknowledgments}
The authors thank Rune Berg for many advice on the neuronal data and Patricia Reynaud-Bouret for fruitful discussions, particularly on the goodness-of-fit part of this work.

\section{Appendix}\label{sec:appendix}

Details on the nonparametric estimators of $b, \sigma^2, a^2$ are developed in \cite{DL2020} and \cite{ADGL2020}. We define here the different estimators used in Section \ref{sec:jumpdiffusion} in Equation (\ref{eq:jumpdiffestim}).
We consider the increments of the observed data:
$$
Y_{k\Delta} := \frac{1}{\Delta} (X_{(k+1)\Delta} - X_{k\Delta})
, \quad 
Z_{k\Delta} := \frac{1}{\Delta} (X_{(k+1)\Delta} - X_{k\Delta})^2.
$$
The smooth truncation function used in the following is 
 $$\varphi(x) =\begin{cases} & 1 \quad |x| <1 \\ & e^{1/3+ 1/(|x|^2-4)}  \\ & 0 \quad |x| \geq 2 \end{cases}.$$
 Let us 
consider the subspaces
  $\mathbb{L}^2(A)$: $\mathcal{S}_m= \text{span}(\varphi_1, \dots, \varphi_{D_m})$
 of dimension $D_m$ 
where $(\varphi_\ell)_\ell$ is an orthonormal basis of $\mathbb{L}^2(A)$.
The subset $A$ is the compact interval of observation (minimum and maximum of observed value).

The nonparametric estimator of $\sigma^2$ is:
$$
\w{\sigma}^2_{\w{m}_\sigma} :=\argmin{m \in \mathcal{M}_n}\argmin{t \in \mathcal{S}_m} \left\{ \frac{1}{n} \sum_{k= 0}^{n - 1} \left(t(X_{k\Delta})- Z_{k\Delta}\varphi\left(\frac{X_{(k+1)\Delta} - X_{k\Delta}}{{\Delta^\beta}}\right)\right)^2 + \pen_\sigma(m)\right\},$$
with 
$\beta \in (0, \frac{1}{2})$ and $\pen_\sigma(m) := \kappa_\sigma \frac{D_m}{n}$.
Then, we define the function $g$:
$$
g(x) := \sigma^2(x) + a^2(x) f(x)
, \quad f(x):= \sum_{j = 1}^M \mathbb{E}[\lambda^{(j)}_{k\Delta} | X_{k\Delta}=x ].
$$
The nonparametric estimator of $g$ is:
$$\w{g}_{\w{m}_g} :=\argmin{m \in \mathcal{M}_n} \argmin{t \in S_m} \left\{ \frac{1}{n} \sum_{k = 0}^{n - 1} ( Y_{k\Delta}-t(X_{k\Delta}))^2+ \pen_g(m)\right\},$$
with $ \pen_g(m) := \kappa_g \frac{D_m}{n\Delta}$.
Finally,
$$\w{a}^2:=\frac{\w{g}_{\w{m}_g}(x)-\w{\sigma}^2_{\w{m}_\sigma}(x)}{\w{f}_{\w{h}}(x)}.
$$
Here
$\w{f}_{\w{h}}$ is the Nadaraya-Watson estimator with a well chosen bandwidth $h$. It is obtained using the estimation of $\lambda$ with $\texttt{ADM4}$.

For the estimation of the drift function, one need the quantity,
$$
 T_{M, k\Delta}= 
 \frac{1}{\Delta} \int_{k\Delta}^{(k+1)\Delta} \sum_{j=1}^M a(X_{s-})d{N}^{(j)}_s,
 $$
 then the nonparametric estimator of $b$ is:
$$ \w{b}_{\w{m}} :=\argmin{m \in \mathcal{M}_n} \argmin{t\in \mathcal{S}_m}\left\{\frac{1}{n} \sum_{k=0}^{n-1} (U_{k\Delta}-t(X_{k\Delta}))^2 +\pen_b(m)\right\}, \quad 
 \pen_b(m): = \kappa_{b}\frac{ D_m}{n\Delta}$$.

%%%%%%%%%%% biblio
%\newpage

\bibliographystyle{plain}
\bibliography{BIB}

%%%%%%%%% FIGURES

%\begin{figure}
%\centering
%\includegraphics[width= 5cm, height=4cm]{images/PSTH_190}
%includegraphics[width= 5cm, height=4cm]{images/PSTH_26Neurons}
%\caption{The PSTH with 190 neurons, then with 26 selected neurons}
%\label{fig:psth}
%\end{figure}
%

%\begin{figure}
%\centering
%\includegraphics[width=5cm, height=4cm]{images/adj_estim_tick_ADM4_seuil01_red_M250_T100_expo}
%\caption{Adjacency matrix for the 12 neurons}
%\label{fig:adjadm}
%\end{figure}
%
%
%
%\begin{figure}
%\centering
%\includegraphics[scale=0.2]{images/histoalltrials}
%\caption{Repartition of the length of the spike trains of the 250 neurons on the concatenated 10 trials}
%\label{fig:histalltrials}
%\end{figure}
%
%\begin{figure}
%\centering
%\includegraphics[scale=0.4]{images/qqplot1trialM26neurons}
%\caption{Goodness-of-fit}
%\label{fig:qqplot}
%\end{figure}

%\begin{figure}
%\centering
%\includegraphics[scale=0.4]{images/plotZtest4M26}
%\caption{Ztest and the theoretical quantile}
%\label{fig:goodnesstest4}
%\end{figure}
%
%
%\begin{figure}
%\centering
%\includegraphics[scale=0.2]{images/qqplotsub3neuronalltrial}
%\caption{Goodness-of-fit plits for residual interarrival times for 3 neurons randomly chosen on the all time interval (concatenated data) qqplotsubsamplealltrials which the subsample technique}
%\label{fig:qqplot}
%\end{figure}
%
%%
%
%
%
%\begin{figure}
%\centering
%\includegraphics[scale=0.2]{images/empiricalIC}
%\caption{Empirical confidence intervals}
%\label{fig:valid}
%\end{figure}

\end{document}